\documentclass[12pt]{article}
\usepackage{amsxtra,amssymb,amsthm,amsmath,latexsym}

\theoremstyle{plain}

\numberwithin{equation}{section}

\def\oH{{\overset{\circ}{H}}}
\def\oH1{{\overset{\circ}{H}\kern-.02in{}^1}}

\def\bee{\begin{equation*}}
\def\eee{\end{equation*}}
\def\be{\begin{equation}}
\def\ee{\end{equation}}

\begin{document}

%\begin{titlepage}
\title{A simple proof of the closed graph theorem}

\author{Alexander G. Ramm \\
Mathematics Department, Kansas State University\\
Manhattan, KS 66506, USA\\
email:     ramm@math.ksu.edu\\
http://www.math.ksu.edu/\,$\widetilde{\ }$\,ramm}

\date{}

\maketitle\thispagestyle{empty}

\begin{abstract}
Assume that $A$ is a closed linear operator defined on all of
 a Hilbert space $H$. Then $A$ is bounded. A new short proof of this classical theorem is given
 on the basis of the uniform boundedness principle. The proof can be easily extended to Banach spaces.
\end{abstract}

 Math subject classification: 47A05

Key words: closed graph theorem; closed linear operator; uniform boundedness principle; new short proof of the closed graph theorem.

\section{INTRODUCTION.}\label{S:1}

We denote by $D(A)$ the domain of definition of $A$, by $A^*$ the adjoint operator,
by $||A||$ the norm of $A$, by $(u,v)$ the inner product in $H$, by $c>0$ various estimation
constants are denoted. Let $A$ be a linear operator in $H$ and $u_n\in D(A)$. Suppose that $u_n\to u$ and $Au_n\to v$.
If the above implies that $u\in D(A)$ and $Au=v$ then the operator $A$ is called closed (see \cite{K}).
 It is well known (see, \cite{DS}, \cite{K}) that any bounded sequence in $H$
contains a convergent subsequence.

The following result is classical:

{\bf Theorem 1.} {\em Let $A$ be a closed linear operator in a Hilbert space $H$,
and $D(A)=H$. Then $A$ is bounded.}

Theorem 1 is known as the closed graph theorem. Its proof
can be found in \cite{DS}, \cite{K}, \cite{Y}, and in many other texts
in functional analysis.

 These proofs are based on the Baire cathegory theorem.
The aim of this note is to give a simple new proof of Theorem 1
using the well-known uniform boundedness principle, which we state as Theorem 2,
and a new result, stated as Theorem 3, which is proved in Section 2.

In \cite{H} a proof of
Theorem 1 is given, which is different from ours.
Our proof of Theorem 1 is not only new but also very short.

 Proofs of Theorem 2  which are not based on Baire's theorem
 can be found in \cite{H}, problem 27,  \cite{He}, \cite{Ho}, \cite{S}.

 {\bf Theorem 2.} {\em If $\sup_n|(Au_n,v)|\le \infty$ for every $v\in H$,
 then $\sup_n||Au_n||\le \infty$}

 We {\em assume} Theorem 2 known.

The new result we use in the proof of Theorem 1 is the following:

{\bf Theorem 3.} {\em If $A$ is a linear closed operator with $D(A)=H$,
then $D(A^*)=H$.}

In Section 2 proofs are given.

 \section{PROOFS.}\label{S:2}

{\bf Proof of Theorem 3}.

If $A$ is a linear closed operator and $D(A)=H$, then $A^*$ exists, is closed and densely defined.
To prove that  $D(A^*)=H$, let $v\in H$ be arbitrary, and $v_n\to v$, $v_n\in D(A^*)$. Let $u\in H$ be arbitrary.
 Then $(Au, v_n)=(u, A^*v_n)$ and
 \begin{equation}\label{eq:1}
\sup_n|(u, A^*v_n)|\le \sup_n ||v_n|| ||Au||\le c(u).
\end{equation}
By Theorem 2 one has $\sup_n||A^*v_n||\le c$. Therefore, a subsequence, denoted again
$A^*v_n$, converges weakly in $H$: $A^*v_n\rightharpoonup v^*$, and $(Au,v)=(u, v^*)$.
Thus, $v\in D(A^*)$, and $D(A^*)=H$ since $v\in H$ was arbitrary.

Theorem 3 is proved. \hfill$\Box$
$$ $$
{\bf Proof of Theorem 1}.  Consider the relation $(Au,v)=(u, v^*)$.
Since $D(A)=H$ and $A$ is closed, Theorem 3 says that $D(A^*)=H$,  the above relation
holds for every $v\in H$, and $v^*=A^*v$.
Suppose that $A$ is unbounded. Then there exists a sequence $u_n$, $||u_n||=1$, such that
\begin{equation}\label{eq:2}
||Au_n||\to \infty.
\end{equation}
  On the other hand, one has:
 \begin{equation}\label{eq:3}
\sup_n|(Au_n, v)|=\sup_n|(u_n,A^*v)|\le \sup_n ||u_n||\cdot||A^*v||=||A^*v||:=c(v).
\end{equation}
 By Theorem 2 one concludes that $\sup_n||Au_n||<c$. This contradicts \eqref{eq:2}.
 Thus, \,\,  one concludes that $||A||<c$. Theorem 1 is proved. \hfill$\Box$

\end{document}